\documentclass[pdflatex,sn-mathphys-num]{sn-jnl}

\usepackage{graphicx}
\usepackage{multirow}
\usepackage{amsmath,amssymb,amsfonts}
\usepackage{amsthm}
\usepackage{mathrsfs}
\usepackage[title]{appendix}
\usepackage{xcolor}
\usepackage{textcomp}
\usepackage{manyfoot}
\usepackage{booktabs}
\usepackage{algorithm}
\usepackage{algorithmicx}
\usepackage{algpseudocode}
\usepackage{listings}

\newcommand{\cH}{\mathcal{H}}
\newcommand{\B}{\mathbb{B}}
\newcommand{\C}{\mathbb{C}}
\newcommand{\R}{\mathbb{R}}
\newcommand{\N}{\mathbb{N}}
\newcommand{\T}{\mathbf{T}}
\newcommand{\W}{\mathcal{W}}
\newcommand{\M}{\mathcal{M}}
\newcommand{\Om}{\Omega}
\newcommand{\cC}{\mathcal{C}}
\newcommand{\cO}{\mathcal{O}}
\newcommand{\cU}{\mathcal{U}}
\newcommand{\Sq}{\mathbb{S}_q}
\newcommand{\co}{\overline{\mathrm{co}}\,}
\newcommand{\cok}{\overline{\mathrm{co}}\,K}
\newcommand{\Fq}{\mathcal{W}_q}
\newcommand{\omq}{\omega_q}
\newcommand{\wq}{w_q}
\newcommand{\csp}{\sigma_p}
\newcommand{\sap}{\sigma_{ap}}

\theoremstyle{thmstyleone}
\newtheorem{theorem}{Theorem}[section]
\newtheorem{proposition}[theorem]{Proposition}
\newtheorem{lemma}[theorem]{Lemma}
\newtheorem{corollary}[theorem]{Corollary}

\theoremstyle{thmstylethree}
\newtheorem{definition}[theorem]{Definition}

\theoremstyle{thmstyletwo}
\newtheorem{remark}[theorem]{Remark}

\begin{document}

\title[Joint $q$-Numerical Ranges]{Joint $q$-Numerical Ranges of Continuous Operator-Valued Functions}

\author*[1]{\fnm{Amin}\sur{Azadi}}\email{aminazadi.edu@gmail.com}

\affil*[1]{\orgname{Independent Research}, \orgaddress{\city{Tehran}, \country{Iran}}}

\articletype{Research Article}

\abstract{This paper develops a unified theory of joint $q$-numerical ranges $\Fq(F)$ for continuous functions $F:\Om \to \B(\cH)$, where $\Om$ is a compact Hausdorff space and $\B(\cH)$ denotes the $C^*$-algebra of bounded linear operators on a complex Hilbert space $\cH$. We prove that $\Fq(F)$ is totally bounded in $\cC(\Om)$, the space of continuous complex-valued functions on $\Om$, and establish the upper semicontinuity of the set-valued map $F\mapsto \Fq(F)$, and show that $\Fq(F)$ is compact when $\cH$ is finite-dimensional. The main result of the paper is a realization theorem. We prove that, for every totally bounded subset $K$ of $\cC(\Om)$ and every $q\in \C$ with $|q|=1$, there exist a Hilbert space $\cH$ and a continuous function $F:\Om\to\B(\cH)$ such that $qK\subset \Fq(F)\subset q\cok$. Consequently, if $K$ is closed and convex, then $\Fq(F)=qK$. As applications, we derive analogous realization results for bounded subsets in $\C^n$ and for totally bounded subsets in arbitrary Banach spaces.}



\maketitle

\section{Introduction and Preliminaries}\label{sec1}

The numerical range is one of the most fundamental and extensively studied concepts in operator theory, with deep connections to spectral theory, stability analysis of dynamical systems, and quantum mechanics. Let $\B(\cH)$ denote the $C^*$-algebra of bounded linear operators on a complex Hilbert space $\cH$. For $T\in \B(\cH)$, the classical numerical range is defined by
\begin{equation}\tag{1.1}
W(T)=\{(Tx,x):x\in \cH,\ \|x\|=1\}.
\end{equation}
The celebrated Toeplitz--Hausdorff theorem \cite[Theorem 1.1.7]{gau2021} asserts that $W(T)$ is always convex, a property that has inspired numerous generalizations and applications across mathematics. For comprehensive treatments of the classical numerical range, we refer the reader to the monographs \cite{gustafson1997,gau2021}. The numerical radius
\[
w(T)=\sup\{|z|:z\in W(T)\}
\]
satisfies the well-known inequalities $\tfrac{1}{2}\|T\|\le w(T)\le \|T\|$.

The concept of the joint numerical range generalizes the classical notion to operator tuples. For an $n$-tuple $\T=(T_1,\ldots,T_n)\in \B(\cH)^n$, the joint numerical range is defined by
\[
\W(\T)=\{\langle \T x,x\rangle:x\in \cH,\ \|x\|=1\},
\]
where
\[
\langle \T x,x\rangle=(\langle T_1x,x\rangle,\ldots,\langle T_nx,x\rangle).
\]
Unlike the classical numerical range, $\W(\T)$ need not be convex when $n\ge 2$; see \cite[Example 3.1]{muller2020}. The geometry and boundary structure of joint numerical ranges have been extensively investigated; see \cite{chan2018,cho1981,li2000} for boundary point characterizations and \cite{gutkin2004} for differential geometric perspectives. The study of joint numerical ranges in more general settings, including semi-Hilbertian spaces and locally convex spaces, has been pursued in \cite{baklouti2018,moore1969,feki2020}.

In a parallel development, Abtahi and Azadi \cite{abtahi2025} introduced and studied the joint numerical range of continuous operator-valued functions. For a compact Hausdorff space $\Om$ and a continuous function $F:\Om\to \B(\cH)$, they defined
\[
\W(F)=\{\langle Fx,x\rangle:x\in \cH,\ \|x\|=1\},
\]
where $\langle Fx,x\rangle(s)=\langle F(s)x,x\rangle$ for all $s\in \Om$. It was shown in \cite[Theorem 2.6]{abtahi2025} that $\W(F)$ is totally bounded in $\cC(\Om)$, the space of all continuous complex functions on $\Om$. They also established the following realization theorem; see \cite[Theorem 2.9]{abtahi2025}.

\begin{theorem}\label{thm:1.1}
Let $K$ be a totally bounded set in $\cC(\Om)$. Then there exist a Hilbert space $\cH$ and an $F\in \cC(\Om,\B(\cH))$ such that
\[
K\subset \W(F)\subset \co K.
\]
Moreover, if $K$ is closed and convex, then $\W(F)=K$.
\end{theorem}

The primary objective of this work is to extend this theorem; the desired extension is established in Theorem \ref{thm:2.16}.

Among the various generalizations of the numerical range, the $q$-numerical range, introduced by Li, Mehta, and Rodman \cite{li1994}, plays a prominent role. Throughout the sequel, let $q\in \C$ satisfy $|q|\le 1$ and $\dim(\cH)\ge 2$. The $q$-numerical range of an operator $T\in \B(\cH)$ is defined as
\[
W_q(T)=\{\langle Tx,y\rangle:(x,y)\in \Sq(\cH)\},
\]
where
\[
\Sq(\cH)=\{(x,y)\in \cH\times \cH:\|x\|=\|y\|=1,\ \langle x,y\rangle=q\}.
\]
It is well-known that $W_q(T)$ is a non-empty subset of $\C$; see \cite{li1994}. A fundamental result of Tsing \cite{tsing1984} establishes that $W_q(T)$ is always convex. The associated $q$-numerical radius
\[
\omq(T)=\sup\{|z|:z\in W_q(T)\}
\]
defines a norm on $\B(\cH)$ for every $q\ne 0$; see \cite[Theorem 3.1]{li1994}. The $q$-numerical ranges of special classes of operators, including normal operators and weighted shift operators, were studied in \cite{chien2005,chien2007}. One has
\begin{equation}\tag{1.1}
\frac{q}{2}\|T\|\le \omq(T)\le \|T\| \qquad (q\in(0,1)).
\end{equation}
Recently, Feki, Patra, Rani, and Taki \cite{feki2026} introduced the joint $q$-numerical range for operator tuples:
\begin{equation}\tag{1.2}
\Fq(\T)=\{\langle \T x,y\rangle:(x,y)\in \Sq(\cH)\},
\end{equation}
where
\[
\langle \T x,y\rangle=(\langle T_1x,y\rangle,\ldots,\langle T_nx,y\rangle).
\]
They investigated its convexity, spectral relations, and numerical radius inequalities. In particular, they established spectral inclusion theorems,
\[
q\csp(\T)\subseteq \Fq(\T),\qquad q\sap(\T)\subseteq \overline{\Fq(\T)},
\]
where $\csp(\T)$ and $\sap(\T)$ denote the joint point spectrum and joint approximate point spectrum, respectively. The study of $q$-numerical radius inequalities has been further developed in \cite{moghaddam2024}. For further results on joint numerical ranges and related topics, we refer the reader to \cite{bonsall1971,feki2020,muller2010,buoni1978,li2009,li2011,li2008,li2009b,muller2018,muller2019,drnovsek2014,popescu2009} and the references therein.

The primary objective of this work is to develop a natural extension of the joint $q$-numerical range of operator tuples to the setting of continuous operator-valued functions.

\subsection*{Notation and conventions}

\begin{enumerate}
\item[(1)] For a compact Hausdorff space $\Om$ and a complex Banach algebra $A$ with identity $\mathbf{1}$, $\|\mathbf{1}\|=1$, let $A^*$ denote the dual space of $A$ and $A_1^*$ be the closed unit ball of $A^*$. We denote by $\cC(\Om,A)$ the Banach algebra of continuous $A$-valued functions $F:\Om\to A$ equipped with the uniform norm
\[
\|F\|_\Om=\sup\{\|F(s)\|:s\in \Om\}.
\]
When $A=\C$, we write $\cC(\Om)$.

\item[(2)] For $f\in \cC(\Om)$ and $a\in A$, we let $fa$ be the $A$-valued function on $\Om$, given by $(fa)(s)=f(s)a$, for all $s\in \Om$.

\item[(3)] For a subset $K$ of $\cC(\Om)$, we denote its uniform closure in $\cC(\Om)$ by $\overline{K}$ and its closed convex hull by $\cok$ (the closure of the convex hull $\mathrm{co}\,K$).

\item[(4)] Let $\M(\Om)$ be the space of complex Radon measures on $\Om$. For $F\in \cC(\Om,A)$ and $\mu\in \M(\Om)$, we define $\mu F=\int_\Om F\,d\mu$. By using \cite[Theorem 3.27]{rudin1991}, $\mu F\in A$ and for every $\phi\in A^*$, we have
\begin{equation}\tag{1.3}
\phi(\mu F)=\mu(\phi\circ F).
\end{equation}
\end{enumerate}

\noindent\textbf{Outline.} This paper develops a unified theory of joint $q$-numerical ranges for $F\in \cC(\Om,\B(\cH))$, extending recent developments in the theory of joint numerical ranges and $q$-numerical ranges. We first establish fundamental properties of $\Fq(F)$, including total boundedness, upper semicontinuity of the set-valued map $F\mapsto \Fq(F)$, and spectral inclusion results and compactness in the finite-dimensional setting. We then prove a realization theorem showing that, for $|q|=1$, every compact convex subset $K$ of $\cC(\Om)$ can be realized as $qK=\Fq(F)$ for a suitable continuous operator-valued function $F$. As immediate applications, we derive analogous realization results for bounded subsets in $\C^n$ and for totally bounded subsets in arbitrary Banach spaces.

\section{Main Results}\label{sec2}

We begin by introducing the following notion.

\begin{definition}\label{def:2.1}
Let $\Om$ be a nonempty set and $q\in \C$ with $|q|\le 1$. For a function $F:\Om\to \B(\cH)$, the joint $q$-numerical range and the joint $q$-numerical radius of $F$ are defined, respectively, by
\begin{equation}\tag{2.1}
\Fq(F)=\{\langle Fx,y\rangle:(x,y)\in \Sq(\cH)\},
\end{equation}
\[
\omq(F)=\sup\{\|f\|_\Om:f\in \Fq(F)\}.
\]
\end{definition}

\begin{remark}\label{rem:2.2}
If $q=1$, then $\W_1(F)=\W(F)$, the joint numerical range defined in \cite{abtahi2025}.
\end{remark}

\begin{remark}\label{rem:2.3}
Every $n$-tuple $\T=(T_1,\ldots,T_n)\in \B^n(\cH)$ can be identified with the function $F_\T:\{1,\ldots,n\}\to \B(\cH)$, defined by $F_\T(i)=T_i$. Under this identification, the joint $q$-numerical range (2.1) of the $n$-tuple $\T$ coincides with the joint $q$-numerical range (1.2) of the function $F_\T$.
\end{remark}

The following two propositions present some basic properties of the joint $q$-numerical range and the joint $q$-numerical radius.

\begin{proposition}\label{prop:2.4}
Let $F,G\in \cC(\Om,\B(\cH))$, $f\in \cC(\Om)$, and $\mu\in \M(\Om)$. Then:
\begin{enumerate}
\item[(1)] $\W_{e^{i\theta}q}(F)=e^{i\theta}\Fq(F)$ for all $\theta\in \R$.
\item[(2)] $\W_{\overline{q}}(F^*)=\overline{\Fq(F)}$, where $F^*(s)=F(s)^*$.
\item[(3)] $\Fq(fF)=\{fg:g\in \Fq(F)\}$.
\item[(4)] $\Fq(F)=\{qf\}$ if and only if $F=fI$, where $I$ is the identity operator.
\item[(5)] $\Fq(F+G)\subset \Fq(F)+\Fq(G)$.
\item[(6)] $\Fq(fI+G)=\{qf+g:g\in \Fq(G)\}$.
\item[(7)] $\Fq(U^*FU)=\Fq(F)$ for all unitary operators $U$.
\item[(8)] $\Fq(\mu F)=\{\mu(f):f\in \Fq(F)\}$.
\end{enumerate}
\end{proposition}

\begin{proof}
We present a proof of (8) only; the remaining statements can be verified directly. For a fixed $(x,y)\in \Sq(\cH)$, define $\phi_{x,y}:\B(\cH)\to \C$ by $\phi_{x,y}(T)=\langle Tx,y\rangle$ for all $T\in \B(\cH)$. Using (1.3), we have
\[
\langle (\mu F)x,y\rangle
=\phi_{x,y}(\mu F)
=\mu(\phi_{x,y}\circ F)
=\mu\langle Fx,y\rangle.
\]
Therefore, $\Fq(\mu F)=\{\mu(f):f\in \Fq(F)\}$.
\end{proof}

\begin{remark}\label{rem:2.5}
If in Proposition \ref{prop:2.4}(8) we take $\mu$ to be a point mass measure at $s$, then for all $s\in \Om$
\begin{equation}\tag{2.3}
W_q(F(s))=\{f(s):f\in \Fq(F)\}.
\end{equation}
\end{remark}

\begin{proposition}\label{prop:2.6}
Let $F\in \cC(\Om,\B(\cH))$.
\begin{enumerate}
\item[(1)] $\omq(F)=\sup_{s\in \Om}\omq(F(s))$.
\item[(2)] The mapping $F\mapsto \omq(F)$ is a seminorm on $\cC(\Om,\B(\cH))$. Moreover, $\omq(\cdot)$ is a norm on $\cC(\Om,\B(\cH))$ if and only if $\wq(\cdot)$ is a norm on $\B(\cH)$.
\item[(3)] $\omq(\cdot)$ defines a norm on $\cC(\Om,\B(\cH))$ for all $q$ with $0<|q|\le 1$.
\item[(4)] Let $q\in(0,1)$. Then
\[
\frac{q}{2}\|F\|_\Om\le \omq(F)\le \|F\|_\Om.
\]
\end{enumerate}
\end{proposition}

\begin{proof}
(1) By (2.3) we have $\wq(F(s))=\sup\{|f(s)|:s\in \Om\}$. Taking the supremum over $s\in \Om$, we get
\[
\sup\{\wq(F(s)):s\in \Om\}
=\sup\{\|f\|_\Om:f\in \Fq(F)\}
=\omq(F).
\]

(2) It is immediate from the definition that $\omq(\cdot)$ is a seminorm. Moreover, a direct computation shows that whenever $\wq(\cdot)$ is a norm, so is $\omq(\cdot)$. It remains to prove the converse implication. Suppose $T\in \B(\cH)$ with $\wq(T)=0$. Let $F:\Om\to \B(\cH)$ be the constant function $F(s)=T$ for all $s\in \Om$. Then $F\in \cC(\Om,\B(\cH))$ and, by part (1), $\omq(F)=\wq(T)=0$. Since $\omq(\cdot)$ is a norm on $\cC(\Om,\B(\cH))$, we have $F=0$, which implies $T=0$. Therefore, $\wq(\cdot)$ is a norm on $\B(\cH)$.

(3) From \cite[Theorem 3.1]{li1994}, $\wq(\cdot)$ defines a norm for all $q$ with $0<|q|\le 1$. The desired result now follows from part (2).

(4) According to (1.1), $\frac{q}{2}\|F(s)\|\le \wq(F(s))\le \|F(s)\|$ for every $s\in \Om$. Now, by taking the supremum over $s\in \Om$ and using part (1), we obtain the desired result.
\end{proof}

The following two theorems investigate the topological properties of the joint $q$-numerical range.

\begin{theorem}\label{thm:2.7}
If $F\in \cC(\Om,\B(\cH))$, then $\Fq(F)$ is totally bounded in $\cC(\Om)$.
\end{theorem}

\begin{proof}
By Proposition 2.1 in \cite{abtahi2025}, $\Fq(F)$ is a bounded set in $\cC(\Om)$. We show equicontinuity. Let $s_0\in \Om$. For any $f=\langle Fx,y\rangle\in \Fq(F)$ and any $s\in \Om$,
\[
|f(s)-f(s_0)|
=|\langle (F(s)-F(s_0))x,y\rangle|
\le \|F(s)-F(s_0)\|.
\]
Since $F$ is continuous at $s_0$, for every $\epsilon>0$ there is a neighbourhood $\cU$ of $s_0$ such that $\|F(s)-F(s_0)\|<\epsilon$ for all $s\in \cU$. Hence $|f(s)-f(s_0)|<\epsilon$ for all $f\in \Fq(F)$. Thus $\Fq(F)$ is equicontinuous at every point. By the Arzela--Ascoli theorem, it is totally bounded.
\end{proof}

The following theorem shows that the set-valued map $F\mapsto \Fq(F)$ is upper semicontinuous in the sense of \cite[Definition 5.15]{bonsall1971}.

\begin{theorem}\label{thm:2.8}
Let $G\in \cC(\Om,\B(\cH))$. For every open set $\cO$ in $\cC(\Om)$ containing $\overline{\Fq(G)}$, there is $\delta>0$ such that if $F\in \cC(\Om,\B(\cH))$ and $\|F-G\|_\Om<\delta$, then $\Fq(F)\subset \cO$.
\end{theorem}

\begin{proof}
By Theorem \ref{thm:2.7}, $\overline{\Fq(G)}$ is a compact subset of $\cC(\Om)$. Hence, there is $\delta>0$ such that if $f\in \cC(\Om)$, $g\in \Fq(G)$, and $\|f-g\|_\Om<\delta$, then $f\in \cO$. Let $F\in \cC(\Om,\B(\cH))$ and $\|F-G\|_\Om<\delta$. If $f\in \Fq(F)$, then $f=\langle Fx,y\rangle$ for some $(x,y)\in \Sq(\cH)$. If $g=\langle Gx,y\rangle$, then
\[
\|f-g\|_\Om
=\|\langle (F-G)x,y\rangle\|_\Om
\le \|F-G\|_\Om<\delta.
\]
Hence $f\in \cO$ and $\Fq(F)\subset \cO$.
\end{proof}

We next define the joint point spectrum and the joint approximate point spectrum associated with $F\in \cC(\Om,\B(\cH))$, and investigate their relationship with the joint $q$-numerical range.

\begin{definition}\label{def:2.9}
Let $F\in \cC(\Om,\B(\cH))$.
\begin{enumerate}
\item[(1)] The joint point spectrum $\csp(F)$ is the set of $f\in \cC(\Om)$ for which there exists a unit vector $x$ of $\cH$ such that $F(s)x=f(s)x$ for all $s\in \Om$.
\item[(2)] The joint approximate point spectrum $\sap(F)$ is the set of $f\in \cC(\Om)$ for which there exists a sequence $(x_k)$ of unit vectors such that
\[
\lim_{k\to\infty}\|F(s)x_k-f(s)x_k\|_\Om=0\qquad \text{for all }s\in \Om.
\]
\end{enumerate}
It is straightforward to verify that $\csp(F)$ and $\sap(F)$ are bounded subsets of $\cC(\Om)$.
\end{definition}

\begin{theorem}\label{thm:2.10}
For $F\in \cC(\Om,\B(\cH))$,
\begin{enumerate}
\item[(1)] $q\csp(F)\subseteq \Fq(F)$;
\item[(2)] $q\sap(F)\subseteq \overline{\Fq(F)}$.
\end{enumerate}
\end{theorem}

\begin{proof}
(1) Let $f\in \csp(F)$. Then there exists a unit vector $x\in \cH$ such that $F(s)x=f(s)x$ for all $s\in \Om$. Choose $z\in \cH$ with $\|z\|=1$ and $\langle x,z\rangle=0$. Set
\[
y=\overline{q}x+\sqrt{1-|q|^2}z.
\]
Then $(x,y)\in \Sq(\cH)$. For each $s\in \Om$,
\[
\langle F(s)x,y\rangle
=\langle f(s)x,y\rangle
=f(s)\langle x,y\rangle
=qf(s).
\]
This implies $qf=\langle Fx,y\rangle$. Therefore $q\csp(F)\subseteq \Fq(F)$.

(2) Let $f\in \sap(F)$. Then there exists a sequence $(x_k)$ of unit vectors of $\cH$ such that
\begin{equation}\tag{2.4}
\|F(s)x_k-f(s)x_k\|\to 0\qquad \text{for all }s\in \Om.
\end{equation}
On the other hand, for each $x_k\in \cH$, there exists $y_k\in \cH$ such that $(x_k,y_k)\in \Sq(\cH)$. In particular, $\langle Fx_k,y_k\rangle\in \Fq(F)$. Then
\[
\|\langle Fx_k,y_k\rangle-qf\|_\Om
=\sup_{s\in \Om}|\langle (F(s)-f(s)I)x_k,y_k\rangle|
\le \sup_{s\in \Om}\|F(s)x_k-f(s)x_k\|.
\]
From (2.4), we get $\langle Fx_k,y_k\rangle\to qf$ in $\cC(\Om)$. Hence $qf\in \overline{\Fq(F)}$.
\end{proof}

It was shown in \cite{feki2026} that $\Fq(\T)$ is a bounded subset of $\C^n$ and is nonempty whenever $\dim(\cH)\ge 2$. In general, however, $\Fq(T)$ need not be closed when $\cH$ is infinite-dimensional. Moreover, Example 2.2 in \cite{feki2026} shows that $\Fq(\T)$ need not be convex. Consequently, $\Fq(F)$ is, in general, neither convex nor closed.

The following theorem provides a sufficient condition for the closedness of $\Fq(F)$.

\begin{lemma}\cite[Proposition 3.7]{feki2026}\label{lem:2.11}
Let $\cH$ be a finite-dimensional complex Hilbert space. Then $\Sq(\cH)$ is compact.
\end{lemma}

\begin{theorem}\label{thm:2.12}
Let $\cH$ be a finite-dimensional complex Hilbert space and $F\in \cC(\Om,\B(\cH))$. Then $\Fq(F)$ is a compact subset of $\cC(\Om)$.
\end{theorem}

\begin{proof}
Let $(f_\alpha)$ be a net in $\Fq(F)$ such that $f_\alpha\to f$ in $\cC(\Om)$. For each $\alpha$, there exists $(x_\alpha,y_\alpha)\in \Sq(\cH)$ such that $f_\alpha=\langle Fx_\alpha,y_\alpha\rangle$. By Lemma \ref{lem:2.11}, $\Sq(\cH)$ is compact. Hence, the net $(x_\alpha,y_\alpha)$ has a convergent subnet $(x_{\alpha_\beta},y_{\alpha_\beta})$ such that $(x_{\alpha_\beta},y_{\alpha_\beta})\to (x,y)$ in $\Sq(\cH)$. Since $F$ is bounded, there exists $M>0$ such that $\|F(s)\|\le M$ for all $s\in \Om$. Then
\begin{align*}
\|\langle Fx_{\alpha_\beta},y_{\alpha_\beta}\rangle-\langle Fx,y\rangle\|_\Om
&=\sup_{s\in \Om}|\langle F(s)x_{\alpha_\beta},y_{\alpha_\beta}\rangle-\langle F(s)x,y\rangle|\\
&\le \sup_{s\in \Om}\bigl(\|F(s)\|(\|x_{\alpha_\beta}-x\|\|y_{\alpha_\beta}\|+\|x\|\|y_{\alpha_\beta}-y\|)\bigr)\\
&\le M(\|x_{\alpha_\beta}-x\|\|y_{\alpha_\beta}\|+\|x\|\|y_{\alpha_\beta}-y\|)\to 0.
\end{align*}
Hence, $\langle Fx_{\alpha_\beta},y_{\alpha_\beta}\rangle\to \langle Fx,y\rangle$ in $\cC(\Om)$. Therefore, $f=\langle Fx,y\rangle\in \Fq(F)$. It follows that $\Fq(F)$ is closed and hence, by Theorem \ref{thm:2.7}, $\Fq(F)$ is a compact subset of $\cC(\Om)$.
\end{proof}

In view of Remark \ref{rem:2.3}, the following is an immediate consequence of Theorem \ref{thm:2.12}.

\begin{corollary}\label{cor:2.13}
Let $\cH$ be a finite-dimensional complex Hilbert space and $\T\in \B^n(\cH)$. Then $\Fq(\T)$ is a compact subset of $\C^n$.
\end{corollary}

Our principal result shows that every compact convex subset of $\cC(\Om)$ can be realized as the joint $q$-numerical range of a suitable continuous operator-valued function. The following two lemmas are needed to establish the main theorem of this section. Since the first lemma is elementary, its proof is omitted for brevity.

\begin{lemma}\label{lem:2.14}
Let $(t_n)_{n=1}^\infty$ be a sequence of complex numbers such that $\sum_{n=1}^\infty |t_n|<\infty$. Suppose that
\[
\Bigl|\sum_{n=1}^\infty t_n\Bigr|=\sum_{n=1}^\infty |t_n|.
\]
Then there exists a complex number $u$ with $|u|=1$ such that $t_n=|t_n|u$ for all $n\in \N$. Equivalently, all non-zero terms have the same argument.
\end{lemma}

\begin{lemma}\label{lem:2.15}
\cite[Lemma 2.8]{abtahi2025} Let $(f_i)$ be a bounded sequence in $\cC(\Om)$. If $(t_i)$ is a sequence in $[0,1]$ such that $1=\sum_{i=1}^\infty t_i$, then $f=\sum_{i=1}^\infty t_i f_i$ belongs to $\overline{\mathrm{co}}\{f_i\}$.
\end{lemma}

We are now ready to state the main theorem, which establishes a $q$-analogue of Theorem \ref{thm:1.1}.

\begin{theorem}\label{thm:2.16}
Let $K$ be a totally bounded set in $\cC(\Om)$, and let $q\in \C$ with $|q|=1$. Then there exist a Hilbert space $\cH$ and an $F\in \cC(\Om,\B(\cH))$ such that
\[
qK\subset \Fq(F)\subset q\cok.
\]
Moreover, if $K$ is closed and convex, then $\Fq(F)=qK$.
\end{theorem}

\begin{proof}
For the proof, we retain the Hilbert space and the function introduced in \cite[Theorem 2.9]{abtahi2025} and apply them directly to the present setting.

Let $\cH=\ell^2(K)$, the space of complex functions $h:K\to \C$ such that $\sum_{f\in K}|h(f)|^2<\infty$. For every $f\in K$, define $e_f:K\to \C$ by
\[
e_f(g)=
\begin{cases}
1,& g=f,\\
0,& g\ne f.
\end{cases}
\]
Then $e_f\in \cH$ and $\|e_f\|_2=1$. For every $s\in \Om$, let $\delta_s$ be the point mass measure at $s$, so that $\delta_s(f)=f(s)$, for every $f\in \cC(\Om)$. Define $F:\Om\to \B(\cH)$ by $F(s)h=\delta_s h$, for all $s\in \Om$ and $h\in \cH$ (here, $\delta_s h$ is the pointwise multiplication of the two functions). It is shown in \cite[Theorem 2.9]{abtahi2025} that $F$ is continuous.

Fix $f\in K$. Let $x=e_f$ and $y=\overline{q}e_f$. Then $\|x\|_2=\|y\|_2=1$ and $\langle x,y\rangle=q\langle e_f,e_f\rangle=q$. Therefore $(x,y)\in \Sq(\cH)$. For each $s\in \Om$,
\[
\langle F(s)x,y\rangle
=\sum_{g\in K}(F(s)e_f)(g)\overline{y(g)}
=q\sum_{g\in K}g(s)|e_f(g)|^2
=qf(s).
\]
This shows that $qf=\langle Fx,y\rangle$. Thus $qK\subset \Fq(F)$.

Conversely, if $f\in \Fq(F)$, then $f=\langle Fx,y\rangle$ for some $(x,y)\in \Sq(\cH)$. Since $x,y\in \ell^2(K)$, there exists a countable set $S=\{f_1,f_2,\ldots\}\subset K$ such that both $x$ and $y$ vanish on $K\setminus S$. Let $t_i=x(f_i)\overline{y(f_i)}$ so that $\sum_{i=1}^\infty t_i=q$. For every $s\in \Om$,
\[
f(s)=\langle F(s)x,y\rangle
=\sum_{i=1}^\infty(\delta_s x)(f_i)\overline{y(f_i)}
=\sum_{i=1}^\infty t_i f_i(s).
\]
By the Cauchy--Schwarz inequality, $\sum_{i=1}^\infty |t_i|\le \|x\|_2\|y\|_2\le 1$. The triangle inequality yields
\[
\Bigl|\sum_{i=1}^\infty t_i\Bigr|\le \sum_{i=1}^\infty |t_i|\le 1.
\]
Since $|q|=1$, we have $\bigl|\sum_{i=1}^\infty t_i\bigr|=\sum_{i=1}^\infty |t_i|=1$. By Lemma \ref{lem:2.14}, there exists $u\in \C$ with $|u|=1$ such that $t_i=|t_i|u$. Hence $q=u$. Put $r_i=|t_i|$, so $\sum_{i=1}^\infty r_i=1$. Thus,
\[
f(s)=\sum_{i=1}^\infty t_i f_i(s)
=q\sum_{i=1}^\infty r_i f_i(s).
\]
This shows that $f=\sum_{i=1}^\infty r_i f_i$. By Lemma \ref{lem:2.15}, $f\in q\cok$. Thus $\Fq(F)\subset q\cok$. Therefore, we obtain $qK\subset \Fq(F)\subset q\cok$. Also, if $K$ is closed and convex, we have $\Fq(F)=qK$.
\end{proof}

In view of Remark \ref{rem:2.3}, we obtain the following result from Theorem \ref{thm:2.16}.

\begin{corollary}\label{cor:2.17}
Let $K$ be a bounded set in $\C^n$, and let $q\in \C$ with $|q|=1$. Then there exist a Hilbert space $\cH$ and a $\T\in \B^n(\cH)$ such that
\[
qK\subset \Fq(\T)\subset q\cok.
\]
Moreover, if $K$ is closed and convex, then $\Fq(\T)=qK$.
\end{corollary}

As a special case of Corollary \ref{cor:2.17}, taking $n=1$, we obtain the following immediate consequence.

\begin{corollary}\label{cor:2.18}
Let $K$ be a bounded set in $\C$, and let $q\in \C$ with $|q|=1$. Then there exist a Hilbert space $\cH$ and a $T\in \B(\cH)$ such that
\[
qK\subset W_q(T)\subset q\co K.
\]
Moreover, if $K$ is closed and convex, then $W_q(T)=qK$.
\end{corollary}

The following theorem is needed to derive the final consequence of Theorem \ref{thm:2.16}.

\begin{theorem}\label{thm:2.19}
\cite[Theorem 1.26]{douglas1998} (Banach) Every Banach space $E$ is isometrically isomorphic to a subspace of $\cC(\Om)$ for some compact Hausdorff space $\Om$.
\end{theorem}

\begin{corollary}\label{cor:2.20}
Let $K$ be a totally bounded set in a Banach space $E$, and let $q\in \C$ with $|q|=1$. Then there exist a Hilbert space $\cH$ and an $F\in \cC(E_1^*,\B(\cH))$ such that
\[
qK\subset \Fq(F)\subset q\cok.
\]
Moreover, if $K$ is closed and convex, then $\Fq(F)=qK$.
\end{corollary}

\begin{proof}
By Theorem \ref{thm:2.19}, every Banach space $E$ is isometrically isomorphic to a subspace of $\cC(\Om)$ for some $\Om=E_1^*$ equipped with the weak$^*$ topology. Define
\[
\Phi:E\to \cC(E_1^*),\qquad (\Phi(x))(\phi)=\phi(x)\quad (x\in E,\ \phi\in E_1^*).
\]
The map $\Phi$ is a linear isometry. Since $K$ is totally bounded in $E$, its image $\Phi(K)$ is totally bounded in $\cC(E_1^*)$. Applying Theorem \ref{thm:2.16} to the set $\Phi(K)$ yields a Hilbert space $\cH$ and a function $F\in \cC(E_1^*,\B(\cH))$ such that
\[
q\Phi(K)\subset \Fq(F)\subset q\overline{\mathrm{co}}\,\Phi(K).
\]
Identifying $K$ with its image $\Phi(K)$ gives the desired inclusions. In addition, if $K$ is closed and convex, then $\Fq(F)=qK$.
\end{proof}

\section*{Conclusion}
This paper introduces the joint $q$-numerical range for continuous operator-valued functions on compact Hausdorff spaces. We establish fundamental properties including total boundedness, upper semicontinuity, and spectral inclusions, and show that compactness follows in finite dimensions. The main result is a realization theorem showing that every totally bounded set in $\cC(\Om)$ can be embedded between the joint $q$-numerical range and its closed convex hull, with equality for closed convex sets. Applications to bounded sets in $\C^n$ and to arbitrary Banach spaces are derived. Potential applications to stability analysis of parameter-dependent operator families and quantum information theory warrant further investigation.

\backmatter

\bmhead{Disclosure statement}
The author has no relevant financial or non-financial interests to disclose.

\bmhead{Funding}
The author declares that no funds, grants, or other support were received during the preparation of this manuscript.

\bmhead{ORCID}
Amin Azadi: \url{https://orcid.org/0009-0001-2118-8319}


\begin{thebibliography}{99}

\bibitem{abtahi2025}
Abtahi M, Azadi A. Joint numerical ranges of continuous operator-valued functions. Acta Sci. Math. (Szeged). 2025. doi:10.1007/s44146-025-00218-1

\bibitem{baklouti2018}
Baklouti H, Feki K, Mahmoud Sid Ahmed OA. Joint numerical ranges of operators in semi-Hilbertian spaces. Linear Algebra Appl. 2018;555:266--284. doi:10.1016/j.laa.2018.06.021

\bibitem{bonsall1971}
Bonsall FF, Duncan J. Numerical Ranges of Operators on Normed Spaces and of Elements of Normed Algebras. London: Cambridge University Press; 1971. doi:10.1017/CBO9781107359895

\bibitem{buoni1978}
Buoni J, Wadhwa B. On joint numerical ranges. Pacific J. Math. 1978;77:303--306. doi:10.2140/pjm.1978.77.303

\bibitem{chan2018}
Chan JT. A note on the boundary of the joint numerical range. Linear Multilinear Algebra. 2018;66:821--826. doi:10.1080/03081087.2017.1326455

\bibitem{chien2005}
Chien MT, Nakazato H. The $q$-numerical ranges of normal operators. Linear Multilinear Algebra. 2005;53(6):393--416. doi:10.1080/03081080500156888

\bibitem{chien2007}
Chien MT, Nakazato H. The $q$-numerical radius of weighted shift operators with periodic weights. Linear Algebra Appl. 2007;422:198--218. doi:10.1016/j.laa.2006.09.017

\bibitem{cho1981}
Cho M, Takaguchi M. Boundary points of joint numerical ranges. Pacific J. Math. 1981;95:27--35. doi:10.2140/pjm.1981.95.27

\bibitem{douglas1998}
Douglas RG. Banach Algebra Techniques in Operator Theory. Grad. Texts in Math., vol. 179. New York: Springer; 1998. doi:10.1007/978-1-4612-1656-8

\bibitem{drnovsek2014}
Drnovsek R, Müller V. On joint numerical radius II. Linear Multilinear Algebra. 2014;62:1197--1204. doi:10.1080/03081087.2013.816303

\bibitem{feki2020}
Feki K. On tuples of commuting operators in positive semidefinite inner product spaces. Linear Algebra Appl. 2020;603:313--328. doi:10.1016/j.laa.2020.06.015

\bibitem{feki2026}
Feki K, Patra A, Rani J, Taki Z. Joint $q$-numerical ranges of operators in Hilbert and semi-Hilbert spaces. Linear Multilinear Algebra. 2026;74(3):294--316. doi:10.1080/03081087.2025.2612152

\bibitem{gau2021}
Gau HL, Wu PY. Numerical Ranges of Hilbert Space Operators. Cambridge: Cambridge University Press; 2021. doi:10.1017/9781108782296

\bibitem{gustafson1997}
Gustafson KE, Rao DKM. Numerical Range: The Field of Values of Linear Operators and Matrices. New York: Universitext, Springer-Verlag; 1997. doi:10.1007/978-1-4613-8498-4

\bibitem{gutkin2004}
Gutkin E, Jonckheere EA, Karow M. Convexity of the joint numerical range: topological and differential geometric viewpoints. Linear Algebra Appl. 2004;376:143--171. doi:10.1016/j.laa.2003.06.011

\bibitem{li1994}
Li CK, Mehta PP, Rodman L. A generalized numerical range: the range of a constrained sesquilinear form. Linear Multilinear Algebra. 1994;37:25--49. doi:10.1080/03081089408818311

\bibitem{li1998}
Li CK, Nakazato H. Some results on the $q$-numerical range. Linear Multilinear Algebra. 1998;43:385--409. doi:10.1080/03081089808818359

\bibitem{li2000}
Li CK, Poon YT. Convexity of the joint numerical range. SIAM J. Matrix Anal. Appl. 2000;21:668--678. doi:10.1137/S0895479898343516

\bibitem{li2009}
Li CK, Poon YT. The joint essential numerical range of operators: convexity and related results. Studia Math. 2009;194:91--104. doi:10.4064/sm194-1-6

\bibitem{li2011}
Li CK, Poon YT. Generalized numerical ranges and quantum error correction. J. Operator Theory. 2011;66(2):335--351.

\bibitem{li2008}
Li CK, Poon YT, Sze NS. Higher rank numerical ranges and low rank perturbations of quantum channels. J. Math. Anal. Appl. 2008;348(2):843--855. doi:10.1016/j.jmaa.2008.08.016

\bibitem{li2009b}
Li CK, Poon YT, Sze NS. Condition for the higher rank numerical range to be non-empty. Linear Multilinear Algebra. 2009;57(4):365--368. doi:10.1080/0308108071786384

\bibitem{lumer1961}
Lumer G. Semi-inner-product spaces. Trans. Amer. Math. Soc. 1961;100:29--43. doi:10.1090/S0002-9947-1961-0133024-2

\bibitem{moghaddam2024}
Moghaddam SF, Mirmostafaee AK, Janfada M. $q$-Numerical radius inequalities for Hilbert space operators. Linear Multilinear Algebra. 2024;72:751--763. doi:10.1080/03081087.2022.2161460

\bibitem{moore1969}
Moore RT. Adjoints, numerical ranges, and spectra of operators on locally convex spaces. Bull. Amer. Math. Soc. 1969;75:85--90.

\bibitem{muller2010}
Müller V. The joint essential numerical range, compact perturbations, and the Olsen problem. Studia Math. 2010;197:275--290. doi:10.4064/sm197-3-5

\bibitem{muller2018}
Müller V, Tomilov Yu. Circles in the spectrum and the geometry of orbits: a numerical ranges approach. J. Funct. Anal. 2018;274(2):433--460. doi:10.1016/j.jfa.2017.10.015

\bibitem{muller2019}
Müller V, Tomilov Yu. Joint numerical ranges and compressions of powers of operators. J. London Math. Soc. 2019;99(1):127--152. doi:10.1112/jlms.12165

\bibitem{muller2020}
Müller V, Tomilov Yu. Joint numerical ranges: recent advances and applications. Concrete Operators. 2020;7(1):133--154. doi:10.1515/conop-2020-0102

\bibitem{popescu2009}
Popescu G. Unitary invariants in multivariable operator theory. Mem. Amer. Math. Soc. 2009;200(941):vi+91.

\bibitem{rudin1991}
Rudin W. Functional Analysis. 2nd ed. New York: McGraw-Hill Inc.; 1991.

\bibitem{tsing1984}
Tsing NK. The constrained bilinear form and the $C$-numerical range. Linear Algebra Appl. 1984;56:195--206. doi:10.1016/0024-3795(84)90125-3

\end{thebibliography}
\end{document}